\documentclass[12pt]{article}
 \usepackage{amsmath}
\usepackage{mathbbol}

\newcommand{\T}{\rm{\scriptscriptstyle T}}
\newcommand{\tr}{ {\rm tr} }
\newcommand{\Exp}[1]{E_{\theta,\Sigma}^{#1}} 
\newtheorem{lemma}{Lemma}
\newtheorem{theorem}{Theorem}

\usepackage{graphicx}
\begin{document}


\begin{center}
{\Large
	{\sc  Scale matrix estimation under data-based loss in high and low dimensions}
}
\bigskip

 Mohamed Anis Haddouche $^{1}$,\, Dominique Fourdrinier $^{2}$ \& Fatiha Mezoued $^{3}$
\bigskip

{\it
$^{1}$ Université de Normandie, INSA Rouen, UNIROUEN, UNIHAVRE, LITIS, avenue de l’Université, BP 8, 76801 Saint-Étienne-du-Rouvray, France.    mohamed.haddouche@insa-rouen.fr et  \'Ecole Nationale Sup\'erieure de Statistique et d\textquoteright
\'Economie Appliqu\'ee (ENSSEA), LAMOPS, Tipaza, Algeria.

$^{2}$ Université de Normandie, UNIROUEN, UNIHAVRE, INSA Rouen, LITIS, avenue de l’Université, BP 12, 76801 Saint-Étienne-du-Rouvray, France.  dominique.fourdrinier@univ-rouen.fr

$^{3}$ \'Ecole Nationale Sup\'erieure de Statistique et d\textquoteright
\'Economie Appliqu\'ee (ENSSEA), LAMOPS, Tipaza, Algeria.\\
famezoued@yahoo.fr
}
\end{center}
\bigskip


{\bf R\'esum\'e.} 
Nous considérons  le problème d'estimation de la matrice d'échelle $\Sigma$ du modèle additif $Y_{p\times n} = M + \mathcal{ E}$, du point de vue de la théorie de la décision. Ici, $p$ représente  le nombre de variables, $n$ le nombre d'observations, $M$ une matrice de paramètres inconnus de rang $q <p$ et  $\mathcal{ E}$  un bruit aléatoire  de distribution à symétrie  elliptique, de matrice de covariance proportionnelle à $I_n \otimes \Sigma$. Ce problème d'estimation est abordé sous une représentation canonique  où la matrice d'observation $Y$ est décomposée en deux matrices, à savoir, $Z_{q\times p}$ qui résume l'information contenue dans $M$ et une matrice $U_{m\times p}$, où $m=n-q$, qui résume  l'information suffisante pour l'estimation de $\Sigma$.
Comme les estimateurs naturels de la forme $  {\hat \Sigma} _a = a \, S $ (où $  S=U^{\T}\,U $ et $ a $ est une constante positive) se comportent mal lorsque  $p > m$ ($S$ n'est pas inversible), nous proposons des estimateurs alternatifs de la forme 
$ {\hat {\Sigma}} _ {a, G} = a \big (  S +  S \,  {S ^ {+} G (Z,S)} \big) $ où ${S^{+}} $ est l'inverse de Moore-Penrose de $ S$ (qui coïncide avec l'inverse $S^{-1}$ lorsque $S$ est inversible). Nous fournissons  des conditions sur la matrice de correction $SS^{+} {G (Z,S)} $ telles que $  {\hat {\Sigma}} _ {a,  {G}} $ améliore  $  {\hat {\Sigma}} _a $ sous le coût basé sur les données $L _S( \Sigma, \hat { \Sigma})
= \tr \big ( S^{+}\Sigma\,({\hat{\Sigma}} \,  {\Sigma} ^ {- 1} -  {I}_ {p} )^{2}\big) $. Nous adoptons une approche unifiée des deux cas où $  S $ est inversible  ($p \leq m$) et $  S $ est non inversible ($p > m$).

{\bf Mots-cl\'es.} Distribution à symétrie elliptique, coût basé sur les données, identité de type Stein-Haff , matrice de covariance, matrice d'échelle. 
\vspace{0.5cm}

\medskip

{\bf Abstract.}
%
We consider  the problem of estimating the scale matrix $\Sigma$ of the additif model  $Y_{p\times n} = M + \mathcal{E}$, under a theoretical decision point of view. Here, $ p $ is the number of variables, $ n$ is the number of observations, $ M $ is a matrix of unknown parameters with rank $q<p$ and   $ \mathcal {E}$ is  a random noise, whose distribution is elliptically symmetric  with covariance matrix  proportional to $ I_n \otimes \Sigma $\,.
We deal  with a canonical form of this model where $Y$ is decomposed in two matrices, namely,  $Z_{q\times p}$ which summarizes the information contained in $ M $, and $ U_{m\times p}$, where $m=n-q$, which summarizes  the sufficient information  to estimate $ \Sigma $.
As the natural estimators of the form ${\hat {\Sigma}}_a=a\,  S$ (where $ S=U^{\T}\,U$ and  $a$ is a positive constant) perform poorly when $p >m$ (S non-invertible), we propose estimators of the form
${\hat{\Sigma}}_{a, G} = a\big(  S+  S \,   {S^{+}\,G(Z,S)}\big)$ where ${S^{+}}$ is the Moore-Penrose inverse of $ S$ (which coincides  with $S^{-1}$ when $S$ is invertible).  We provide conditions on the correction matrix $SS^{+}{G(Z,S)}$ such that ${\hat {\Sigma}}_{a,  G}$ improves over ${\hat {\Sigma}}_a$ under the data-based loss $L _S( \Sigma, \hat { \Sigma}) =\tr \big ( S^{+}\Sigma\,({\hat{\Sigma}} \,  {\Sigma} ^ {- 1} -  {I}_ {p} )^ {2}\big) $.
We adopt a unified approach of the two cases where $ S$  is invertible ($p \leq m$) and $ S$ is non-invertible ($p>m$).


\bigskip\bigskip


\section{Introduction}

Consider  the following additive model 
\begin{align}\label{additive.model}
{Y}={M}+  {\mathcal{E}}  ,\hspace{1.cm} {\mathcal{E}}  \sim  {ES}  ({0}_{np},I_n\otimes \Sigma), 
\end{align}
where ${Y}$ is an observed $n\times p$ matrix, ${M}$ denotes an  $n\times p$ matrix of  unknown parameters and ${\mathcal{E}}$ is an $n\times p$ elliptically symmetric distributed noise with unknown covariance matrix proportional to $I_n\otimes \Sigma$, where $\Sigma $ is an unknown $p\times p$   invertible scale matrix and $I_n$ is the $n$-dimensional identity matrix.
Note that, the class of  elliptically symmetric distributions encompasses a large number of important distributions such as Gaussian, Cauchy, exponential, Student and   Weibull distributions. Our main assumption is that $ M$ is of low-rank, that is,
\begin{align}\label{low-rank}
{\rm rank}( M)= q < p 
\end{align}
%

Note that Model (\ref{additive.model}) is a common alternative representation of the multivariate low-rank regression model $Y = X\,{\beta} + \mathcal{E}$, where $X$ is an $n\times q$ matrix of known constants of rank $q < p$ and ${\beta}$ is an $q\times p$ matrix of unknown parameters. 
In the Gaussian setting, Model (\ref{additive.model}) arises in many fields that require to estimate $M$ as in signal processing, image processing, collaborative filtering. Thus, it has been considered by various authors such as  Cand\`es and Recht (2009), Ji et al. (2010) and Cand\`es et al. (2013).
Recently,  Canu and Fourdrinier (2017) introduced the extended elliptical setting in Model (\ref{additive.model}).
It is worth noting that many estimation procedures of $M$ rely on an accurate estimation of the scale matrix $\Sigma$, which is the aim of this paper.

Thanks to the low-rank assumption in  (\ref{low-rank}), there exists a $n\times n$ orthogonal matrix $Q=(Q_1 Q_2)$, with $Q_2^{\T}\,M = 0$, so that  the canonical form of Model (\ref{additive.model}) is given by 
\begin{align}\label{canonical.model}
Q^{\T}\,Y=
\begin{pmatrix} 
{Q_1^{\T}} \\ {Q_2^{\T}} 
\end{pmatrix}
\,Y
= 
\begin{pmatrix} 
{Z} \\ {U} 
\end{pmatrix}
= 
\begin{pmatrix} 
{\theta} \\ {0} 
\end{pmatrix}
+ Q^{T }\,{\mathcal{E}},
\end{align}	
where $ Z$ and $ U$ are respectively $q\times p$ and $m\times p$ matrices (\textit{cf.} Fourdrinier and Canu~(2017)  for more details). Note that, the canonical form $(\ref{canonical.model})$ separates information about the mean structure $Z$ and the information concerning the scale $U$, since $ S=U^{T}\,U$   summarizes  the   information to estimate $ \Sigma $. 
Now, we restrict our attention to the setting where   the joint density of ${Z}$ and ${U}$ is of the form
\begin{align}\label{density}
({z},{u}) 
&\mapsto
|\Sigma|^{-n/2} \,f \! \left[\,\tr \{({z} 
- {\theta})\Sigma^{-1}({z}-
{\theta})^{\T }\} + \tr \{\Sigma^{-1}\,{u}^{\!\T }
{u}\}\,\right]\,,
\end{align}
for some function $f$. 

In the following, $\Exp{}$ will denote the expectation with respect to the density  (\ref{density}) and $E^{\ast}_{\theta,\Sigma}$ the expectation with respect to the density
\begin{align*}
({z},{u}) \mapsto \frac{1}{K^{*}}|\Sigma|^{-n/2}\,F^{*}
\! \left[\,\tr\{\,({z} - 	
{\theta})\Sigma^{-1}({z}-{\theta})^{\top}\} + \tr\{\Sigma^{-1}\,u^{\T}\,u\}\,\right]\,,
\end{align*}  
where
${F^{*}(t) =\frac{1}{2} \, \int^{\infty}_{t} f(\nu) \,d\nu}$  and the normalizing constant $K^{*}$ is assumed to be finite. Note that, in the setting of a multivariate normal distribution, since $F^{*}=f$,
these two expectations coincide.

As mentioned by James and Stein (1961),  the natural estimators of the form  ${\hat {\Sigma}}_a=a\,  S$ (where  $a$ is a positive constant) perform poorly. Therefore, we consider  alternative estimators   of the form $\hat{\Sigma}_{{a},G}= a\,(S + SS^{+}G(Z,S))$ and we derive dominance results under the data-based loss function
\begin{align}\label{loss}
L_{S}(\hat{\Sigma},\Sigma) = \tr\big(S^{+}\Sigma\,\big({\hat\Sigma}\Sigma^{-1}- I_{p}\big)^{2}\big),
\end{align}
and its associated risk
\begin{align}\label{risk}
R(\hat{\Sigma},\Sigma)
=\Exp{} [L_S(\hat{\Sigma},\Sigma)]\,,
\end{align}
where $\hat{\Sigma}$ is an estimator of $\Sigma$ and $S^{+}$ is the Moore-Penrose inverse of $S$. It is worth  noticing that this type of loss function, called data-based since it involves $ S ^ {+} $, was introduced by Efron and Morris (1976). Since then, it was considered by various authors, in a Gaussian setting by Kubokawa and Srivastava (2008) and Tsukuma and Kubokawa (2015), and in a spherical setting by Fourdrinier and Strawderman (2015).

%

The two main features of our approach is that we consider the general elliptically symmetric distribution context and we unify the two cases where $S$ is non-invertible ($p>m$) and $S$ is invertible ($p \leq m$). The primary decision-theoretic results are presented in Section \ref{Improved estimators}. More precisely, we derive a sufficient condition on the correction matrix function $SS^{+}G(Z,S)$ for which  $\hat{\Sigma}_{{a},G}$ improves on $\hat{ \Sigma}_{a}$ under the data-based loss in  (\ref{loss}). In Section \ref{Simulation}, we provide  numerical results through simulations.


 

\section{Main result}\label{Improved estimators}
Among the  usual estimators $\hat{\Sigma}_a=a\,S$, there exists $a_o >0$ such that $\hat{\Sigma}_{a_o}$ is optimal (that is, the risk of $\hat{\Sigma}_{a_o}$ is less than or equal to the risk of $\hat{\Sigma}_{a}$, for any $a>0$); this is
\begin{align*}
a_o = \frac{1}{K^{*}(p \vee m)}\,,
\end{align*}
where $p\vee m = \max\{p,m\}$ (\textit{cf.} Haddouche (2019) for a proof).
The improvement over the class of $a\,S$'s will be shown through the improvement of
\begin{align}\label{alternative.esti}
\hat{\Sigma}_{{a_o},G}= a_o\,(S + SS^{+}G(Z,S))\,,
\end{align}
over $\hat{ \Sigma}_{a_o}=a_o\,S$,
where $$G(Z,S)=\frac{t}{\tr(S^{+})}\,SS^{+}$$ and $t$ is a positive constant. Note that the choice of this specific form of $G(Z,S)$ is motivated by the estimator considered by Konno (2009) in the normal case.
We give sufficient conditions on the corrected factor $SS^{+}G({Z},S)$, that is on the constant $t$, such that the risk difference
\begin{align*}
\Delta(G) = R(\hat{\Sigma}_{{a_o},G}, \Sigma) - R(\hat{\Sigma}_{a_o} , \Sigma)
\end{align*}
between $\hat{\Sigma}_{{a_o},G}$ and $\hat{\Sigma}_{a_o}$  is non-positive. Of course, $\Delta(G) \leq 0$ makes only sense if and only if $R(\hat{\Sigma}_{{a_o},G}, \Sigma) < \infty$. It is shown in Haddouche (2019) that this occurs as soon as the expectations $\Exp{}\left[\|S^{+}\,G \|^{2}_F\right]
$, $\Exp{}\left[\|\Sigma^{-1}SS^{+}G \|^{2}_F\right]
$, $\Exp{} \left[\tr (\Sigma\,S^{+})\right]$ and  $\Exp{} \left[\tr (\Sigma^{-1}\,S)\right]$ are finite. In that case,
\begin{align}\label{Delta.SG.1.S}
\Delta(G)
&=
a_{o}^2\,K^{*}\,\Exp{}\big[\tr\left(\Sigma^{-1}SS^{+}G\left\{I_p+ S^{+}G + SS^{+}	\right\}\right)\big]
-2\,a_{o}\,\Exp{}\big[\tr\left(S^{+}\,G\right)\big]\,.
\end{align}
%
%

The dependence of the risk difference in (\ref{Delta.SG.1.S}) on the unknown parameter $\Sigma^{-1}$ is
problematic. As a remedy, we apply the  Stein-Haff type identity in the framework of elliptically symmetric distribution given in Fourdrinier Haddouche and Mezoued~(2019).
\begin{lemma}
	Let $G(z,s)$ be a $p \times p$ matrix function such that, for any fixed $z$, $G(z,s)$ is weakly differentiable with respect to $s$. Assume that $\Exp{} \left[\big| \tr \left(\Sigma^{-1} S \,S^{+}\,G\right) \big| \right] < \infty$. Then we have 
	\begin{align*}
	\Exp{}\left[\tr\left(\,\Sigma^{-1}\,SS^{+}\,G\right)\right] 
	=
	K^{*}\,\Exp{*}\left[\tr \left(
	2\,SS^{+}\,{\cal D}_s
	\{SS^{+}G\}^{\T}\, + (m-(p\wedge m)-1)\,S^{+}\,G\right)\,\right].
	\end{align*}
\end{lemma}

Thanks to this identity, sufficient conditions for improvement of
$\hat{\Sigma}_{a_o,G}$ over $\hat{\Sigma}_{a_o}$, are given in the following
theorem  (\textit{cf.} Haddouche (2019) for a proof) through an upper bound of the risk difference in (\ref{Delta.SG.1.S}).

\begin{theorem}
	Consider a density of the form (\ref {density}). Let
	\begin{align}\label{Est.Konno.New}
	\hat{\Sigma}_{a_o,G} =  a_o\left( S + \frac{t}{\tr(S^{+})}\, SS^{+} \right) 
	\end{align}
	where $t$ is a positive constant. Then $\hat{\Sigma}_{a_o,G}$ improves over  $\hat{\Sigma}_{a_o}$ as soon as 
	\begin{align*}
	0 \leq t
	\leq \frac{2\,((p\wedge m) -1 ) }{(p\vee m) - (p \wedge m) + 1}\,.
	\end{align*}
where $p \wedge m = \min\{p,m\}$\,.
\end{theorem}


\section{Numerical study}\label{Simulation}
We deal here with the non-invertible case ($p>m$) for  a Gaussian distribution ($K^{*}=1$) where the scale matrix have an autoregressive structure of the form $(\Sigma)_{ij} = 0.9^{|i-j|}$. Note that simulation on the Student distributions are under study. We evaluate numerically the performance of the alternative estimator $\hat{\Sigma}_{a_o,G}$ in (\ref{Est.Konno.New}) where ${a_o}=1/p$ and $t= 2\,(m-1)/(p-m+1)$, through  the  percentage relative improvement in average loss PRIAL of $\hat{\Sigma}_{a_o,G}$ over $\hat{\Sigma}_{a_o}$ defined as 
\begin{align*}
\text{PRIAL}(\hat{\Sigma}_{a_o,G})= \frac{\text {average loss of }\,\,\hat{ \Sigma}_{a_o}- \text {average loss of}\,\, \hat{\Sigma}_{a_o,G} }{\text {average loss of} \,\,\hat{ \Sigma}_{a_o}}\, \times 100\,,
\end{align*}
which is reported in the following table. 
\vspace{.5cm}

	\begin{tabular*}
		{1\textwidth}{@{\extracolsep{\fill}}  c  c  c  }
		\hline $p$& $m$ &PRIAL\,(\%)   \\
		\hline 
		20 & 4 	 &15.00  \\ 		
		20 & 8 	  &18.56  \\ 	
		20 & 12   &25.56  \\ 
		20 & 16   &47.034  \\  
		100 & 20    &3.39  \\ 
		100 & 40    &4.19  \\  
		100 & 60    &5.76  \\ 
		100 & 80    &10.42  \\ 
		\hline 
	\end{tabular*} 
		 \begin{center}
Results of 1000 Monte Carlo simulation for $(\Sigma)_{ij} = 0.9^{|i-j|}$.
			\end{center}  
For $p=20$ and $p=100$, the PRIAL increases with the values of $m$. Note that, when $p=20$ and $m=16$ the PRIAL is close de $50\,\%$. Note that the data-based Loss is much more discriminant then the usual quadratic loss for which the PRIAL is lower.

\newpage
\section*{Bibliographie}

\noindent Cand{\`{e}}s, E. J. and  Sing{-}Long, C. A. and  Trzasko, J. D. (2013).  Unbiased Risk Estimates for Singular Value Thresholding and Spectral Estimators.\textit{{IEEE} Transactions on signal processing} 61:4643-4657 

\noindent Cand{\`e}s, E. J. and Recht, B. (2009). Exact Matrix Completion via Convex Optimization, \textit{Foundations of Computational Mathematics}, 9(6): 717. 

\noindent Efron, B. and Morris, C. (1976). Multivariate empirical Bayes and estimation of covariance matrices. \textit{Annals of Statistics.} 4(1):22-32 

\noindent Fourdrinier, D. and Canu, S. (2017). Unbiased risk estimates for matrix estimation in the elliptical case, \textit{Journal of Multivariate Analysis}, 158, pp. 60-72

\noindent Foudrinier, D. Haddouche, M. A. and Mezoued, F. (2019).  Scale matrix estimation of an elliptically symmetric distribution in high and low dimensions, \textit{ Université de Rouen Normandie et ENSSEA Tipaza,Technical report} 

\noindent Fourdrinier, D. and Strawderman, W.E. (2015). Robust minimax {S}tein estimation under invariant data-based loss for spherically and elliptically symmetric distributions. \textit{Metrika.} (4)78:461-484 

\noindent Haddouche, M.A. (2019). Scale matrix estimation. \textit{Ph.D disertation, Normandie Universit{\'e}}. Chapter 4 

\noindent James, W. and Stein, C. (1961). Estimation with Quadratic Loss. \textit{Proceedings of the Fourth Berkeley Symposium on Mathematical Statistics and Probability, Volume 1: Contributions to the Theory of Statistics.} 361--379

\noindent  Ji, H. and Liu, C. and Shen, Z. and Xu, Y. (2010). Robust video denoising using low rank matrix completion  \textit{2010 IEEE Computer Society Conference on Computer Vision and Pattern Recognition.} 100:2237-2253

\noindent Kubokawa, T. and Srivastava, M.S. (2008). Estimation of the precision matrix of a singular {W}ishart distribution and its application in high-dimensional data. \textit{Journal of Multivariate Analysis.} (9)99:1906-1928 	

\noindent Konno, Y. (2009). Shrinkage estimators for large covariance matrices in multivariate real and complex normal distributions under an invariant quadratic los. \textit{Journal of Multivariate Analysis} 100:2237-2253 

\noindent Tsukuma, H. and Kubokawa, T. (2015). A unified approach to estimating a normal mean matrix in high and low dimensions. \textit{Journal of Multivariate Analysis} 139:312-328

\end{document}